\documentclass[12pt,a4paper,twocolumn]{article}
\usepackage{amsmath,amssymb,amsfonts,amsthm}
\usepackage{mathrsfs} 
\usepackage{graphicx}
\newcommand{\G}{\mathscr G}
\newcommand{\R}[1]{\mathscr{R}(#1)}
\newcommand{\N}[1]{\mathscr{N}(#1)}
\newcommand{\ov}[1]{\mathscr D(#1)}
\newcommand{\mof}{\mathscr{F}}

\newcommand{\oD}{\mathcal D}
\newcommand{\oB}{\mathcal B}
\newcommand{\oH}{\mathcal H}

\newcommand{\dt}{\mathrm dt}

\newcommand{\f}[1][\cdot]{f(#1)}
\newcommand{\z}{z(\cdot)}

\newcommand{\hp}[1][\cdot]{\hat{p}(#1)}
\newcommand{\hpt}{\hp[t]}
\newcommand{\x}[1][\cdot]{x(#1)}
\newcommand{\xt}{x(t)}
\newcommand{\hx}[1][\cdot]{\hat{x}(#1)}
\renewcommand{\u}{u(\cdot)}
\newcommand{\hu}[1][\cdot]{\hat{u}(#1)}
\newcommand{\hut}{\hu[t]}
\newcommand{\el}{\ell(\cdot)}

\newcommand{\Lr}{\mathbb{L}_2^r}
\newcommand{\Ln}{\mathbb{L}_2^n}
\newcommand{\Lm}{\mathbb{L}_2^m}

\newcommand{\LN}{\mathbb{L}_2([0,\omega],\Rn)}

\newcommand{\Rn}{\mathbb{R}^n}

\theoremstyle{definition}
\newtheorem{dfn}{Definition}
\theoremstyle{theorem}
\newtheorem{thm}{Theorema}

\newtheorem{col}{Corollary}

\begin{document}
\begin{center}
\textbf{MINIMAX STATE OBSERVATION IN LINEAR ONE DIMENSIONAL 2-POINT BOUNDARY VALUE PROBLEMS}\\[5pt]
  Serhiy Zhuk, Serhiy Demidenko, Alexander Nakonechniy\\
        Faculty of cybernetics\\Taras Shevchenko Kyiv
    National University, Ukraine\\e-mail: beetle@unicyb.kiev.ua
\end{center}

\textbf{Abstract.}  In this paper we study observation problem for linear 2-point BVP $\oD\x=\oB\f$ assuming that information about system input $\f$ and random noise $\eta$ in system state observation model $y(\cdot)=\oH\x+\eta$ is incomplete ( $\f$ and $M\eta\eta'$ are some arbitrary elements of given sets). A criterion of guaranteed (minimax) estimation error finiteness is proposed. Representations of minimax estimations are obtained in terms of 2-point BVP solutions. It is proved that in general case we can only estimate a projection of system state onto some linear manifold $\mof$. In particular, $\mof=\Ln$ if  $\mathrm{dim}\mathscr{N}\bigl(\begin{smallmatrix}\oD\\\oH\end{smallmatrix}\bigr)=0$. Also we propose a procedure which decides if given linear functional belongs to $\mof$.
\normalsize \small{
\section*{\textbf{Problem statement}}
}

\normalsize
Let $t\mapsto\xt$  -- totally continuous vector-function from space of square summable $n$-vector-functions $\Ln:=\LN$ -- be a solution of BVP
\begin{equation}
  \label{eq:kr}
  \dot x(t)-A(t)x(t)=B(t)f(t),x(0)=x(\omega),
\end{equation}
where $t\mapsto A(t)$($t\mapsto B(t)$) -- $n\times n$($n\times r$)-matrix-valued continuous function, $\omega<+\infty$, $\f\in\Lr$. 

We suppose that a realization of $m$-vector function $t\mapsto y(t)$ is observed at $[0,\omega]$ 
\begin{equation}
  \label{eq:y}
  y(t)=H(t)x(t)+\eta(t),
\end{equation}
where $t\mapsto x(t)$ is one of the possible solutions of \eqref{eq:kr} for some $\f\in\G$, $t\mapsto H(t)$ -- $m\times n$-matrix-valued continuous function, $t\mapsto\eta(t)$ -- realization of mean-square continuous random process with zero expectation and uncertain correlation function $(t,s)\mapsto R_\eta(t,s)\in\G_2$. Let $$
\G:=\{\f:\int_0^\omega(f(t),f(t))\dt\le1\},$$$$
\G_2:=\{R_\eta:\int_0^\omega\mathrm{sp}R_\eta(t,t)\dt\le1\}
$$ and consider linear functional $$
\ell(x):=\int_0^\omega(\ell(t),x(t))\dt,\quad\el\in\Ln,
$$ defined on the~\eqref{eq:kr} solutions domain. \emph{We will be looking for $\ell(x)$ estimation in terms of} $$
u(y):=\int_0^\omega(u(t),y(t))\dt,\quad\u\in U_\ell\subset\Lm
$$ For each $\u$ we associate \emph{guaranteed estimation error}\footnote{Linear mapping $\oD$ is defined by the rule $
\oD x=\dot x-Ax, x\in\ov\oD,
$ where $\ov\oD$ is set of totally continuous vector-functions $t\mapsto x(t)$ satisfying $
\int_0^\omega|\dot x(t)|^2_n<+\infty,\int_0^\omega\dot x(t)\dt=0
$, $x\mapsto Ax$ multiplies $\x$ by $t\mapsto A(t)$.} 
 $$
 \sigma(u):=\sup_{\x\in\ov\oD,\oD \x\in\G,R_\eta\in\G_2}\{M[\ell(x)-u(y)]^2\}
$$
\begin{dfn}
  Function $\hu\in U_\ell$ is called \emph{minimax mean-square estimation} if it satisfies 
  \begin{equation}\label{err}
     \sigma(\hat{u})\le
     \sigma(u),\quad \u\in U_\ell
  \end{equation}
  Term $$
  \hat{\sigma}:=\inf_{u\in U_l}\sigma(u)
$$ is called minimax mean-square error.
\end{dfn}
\begin{thm}\label{t:3}
Boundary value problem
  \begin{equation}
  \label{eq:hphz}
  \begin{split}
    &\dot z(t)=-A'(t)z(t)+H'(t)H(t)p(t)-\ell(t),\\
    &\dot p(t)=A(t)p(t) + B(t)B'(t)z(t),\\
    &z(0)=z(\omega),p(0)=p(\omega)\\
  \end{split}
\end{equation}
has non-empty solutions domain iff $$
P h(\omega)\perp\N{W(0,\omega)},
$$
where $P:=[E-(E-\Phi(\omega,0))(E-\Phi(\omega,0))^+]$, $\Phi$ -- fundamental solution of $\dot z(t)=-A'(t)z(t)$,
$$
W(0,\omega):=\int_0^\omega P\Phi(\omega,s)H'(s)H(s)\Phi'(\omega,s)P\mathrm{ds},
$$
$h(\cdot)$ is a solution of $$
\dot h(t)=-A'(t)h(t)+\ell(t),h(0)=0
$$
\end{thm}
Let's illustrate theorem ~\ref{t:3}. Set
$$
A(t)\equiv\begin{pmatrix}1&0\\1&0\end{pmatrix},
B(t)\equiv\begin{pmatrix}1&0\\0&1\end{pmatrix},
H(t)\equiv\begin{pmatrix}1&0\\0&0\end{pmatrix}
$$
Fundamental solution $t\mapsto F(t)$ of~\eqref{eq:kr} (and fundamental solution $t\mapsto G(t)$ of adjoint BVP)
$$
F(t)\equiv\begin{pmatrix}e^{t}&0\\-1+e^{t}&1\end{pmatrix},
G(t)\equiv\begin{pmatrix}e^{-t}&e^{-t}-1\\0&1\end{pmatrix}
$$
than $\N\oD=\{(0,1)\}$ and $\oH\N\oD=\{0\}$.
Let $\el=l_1(\cdot)=\bigl[\begin{smallmatrix}\sin(t)\\1\end{smallmatrix}\bigr]$. Than 
$$h(t)=\bigl[\begin{smallmatrix}-\frac{1}{2}e^{-t}(1-2e^t+2e^t t+e^t\cos(t)-e^t\sin(t))\\t\end{smallmatrix}\bigr]$$
and
$$
P=\begin{pmatrix}0&0\\0&1\end{pmatrix},
W(2\pi,0)\equiv\begin{pmatrix}0&0\\0&0\end{pmatrix}
$$
As far as $W(2\pi,0)$ is a zero matrix, than according to theorem~\ref{t:3} $\el\in\mof$ if and only if $Ph(2\pi)=0$. But for chosen $l_1(\cdot)$$$
h(2\pi)=\bigl[\begin{smallmatrix}\frac{1}{2}-\frac{e^{-2\pi}}2-2\pi\\2\pi\end{smallmatrix}\bigr]\Rightarrow
Ph(2\pi)=\bigl[\begin{smallmatrix}0\\2\pi\end{smallmatrix}\bigr]
$$ 
Let $\ell(t):=l_2(t)=(\sin(t),\cos(t))$. Than $$h(t)=(0,\sin(t))\Rightarrow Ph(2\pi)=(0,0)$$
It's easy to see that \eqref{eq:hphz} solution's domain is empty for $(0,l_1(\cdot))$. Really, null-space of adjoint BVP is $N=\{(0,0,0,1)\}$ and $(0,l_1(\cdot))$ is not orthogonal to $N$ while $(0,l_2(\cdot))\perp N$.

Let's denote by $\mof$ set of all $\el\in\Ln$ satisfying condition of the theorem~\ref{t:3}. In the next theorem we state that minimax error is finite iff  $\el\in\mof$ and in that case unique minimax estimation $\hu$ exists. 
\begin{thm}\label{t:1}
Minimax mean-square error
$$
\hat{\sigma}=
\begin{cases}
  +\infty,&\el\notin\mof,\\
  \int_0^\omega(\ell(t),\hpt)_n\dt
\end{cases}
$$ If $\el\in\mof$ than unique minimax estimation $\hu$ exists and $$
\hut=H(t)\hpt,
$$ where $\hp$ is one of the \eqref{eq:hphz} solutions.
\end{thm}
\begin{col}\label{t1c1}
For given $y(\cdot)\in\Lm$ minimax estimation $\hu$ can be represented as $$            
\int_0^\omega(\hut,y(t))\dt=\int_0^\omega(\ell(t),\hat{x}(t))\dt,
$$ where $\hx$ is any solution of
\begin{equation}
  \label{eq:hxhx}
  \begin{split}
    &\dot p(t)=-A'(t)p(t) - H'(t)(y(t)-H(t)x(t)),\\
    &\dot x(t)=A(t)x(t)+B(t)B'(t)p(t),\\
    &p(0)=p(\omega),x(0)=x(\omega)
  \end{split}
\end{equation}
\end{col}
\begin{col}
If system of functions\footnote{$\oH\psi_k(t)=H(t)\psi_k(t)$, $\psi_k(\cdot)$ are linearly independent solutions of the homogeneous BVP~\eqref{eq:kr}. }  $\{\oH\psi_k(\cdot)\}$ is linear independent, than for all $\el\in\Ln$ minimax estimation is represented in terms of theorem~\ref{t:1} or previous corollary.
\end{col}
\begin{col}
If $L$ is linear Noether closed mapping in $\Ln$, $\oH,\oB$ are bounded linear mappings in $\Ln$ than $$
(0,\ell)\in\R{
  \begin{smallmatrix}
  -L&&\oB\oB'\\\oH'\oH&&L'
  \end{smallmatrix}}\Leftrightarrow\el=L'z+\oH'\u
$$ 
for some $\z,\u\in\Ln$.
\end{col}
%
\textbf{Example 1.}
We will apply corollary~\ref{t1c1} to linear oscillator's state estimation problem
$$
A(t)\equiv\begin{pmatrix}0&-1\\1&0\end{pmatrix},
B(t)\equiv\begin{pmatrix}1&0\\0&1\end{pmatrix},
$$$$
H(t)\equiv\begin{pmatrix}\frac{\cos t}{20}&\frac{\sin t}{20}\\\frac{\cos t}2&
\frac{\sin t}2
\end{pmatrix}
$$ It's easy to see that $$
\N\oD=\{\{\cos(t),-\sin(t)\},\{\sin(t),\cos(t)\}\},
$$ hence $$
\oH\N\oD=\{\{0,0\},\{\frac 1{20},\frac 12\}\}
$$
Let $f(t)=\begin{pmatrix}\frac{\cos(t)}{\pi}\\\frac{\sin(t)}{\pi}\end{pmatrix}$ and suppose $$x(t)=\begin{smallmatrix}
\cos(t)/2 +\sin(t) + t\sin(t)/\pi\\
\cos(t) + t\ cos(t)/\pi - \sin(t)/2)\end{smallmatrix}
$$ is observed while noise $g(t)=\begin{pmatrix}0.1\sin(t)\\0.1\sin(t)\end{pmatrix}$. Than output $y(t)=((0.05 + 0.0159155 t + 0.1\sin(t),   0.5 + 0.159155 t + 0.1\sin(t))$, so we do not have any info about component from $\oD$ kernel $(\cos(t)/2,- \sin(t)/2)$ included in $x(t)$.\\
Let's find $\hx$ from~\eqref{eq:hxhx}. We obtain $$
\|\x-\hx\|_2\simeq 1.85877
$$ and ($\x$ -- solid line, $\hx$ -- dashed line)
\begin{center}
\includegraphics{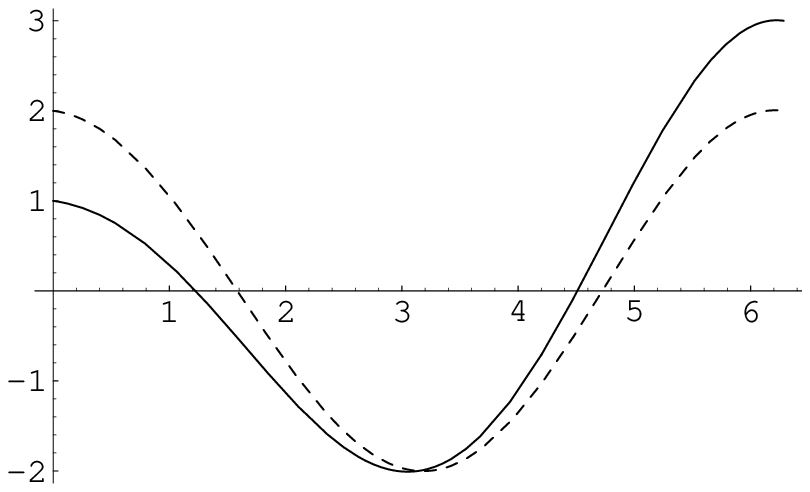}\\
\includegraphics{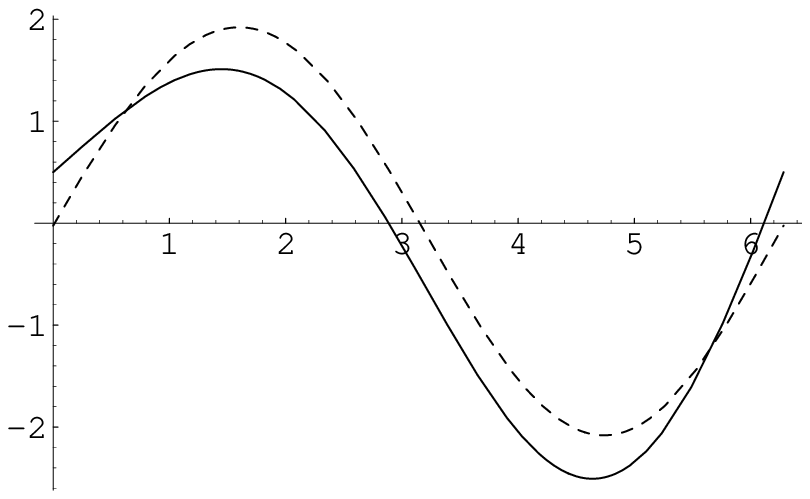} 
\end{center}
According to theorem~\ref{t:1} in general case we can only estimate a projection of~\eqref{eq:kr} state onto linear manifold $\mof$. In particular, if $\N\oH\cap\N\oD={0}$, than $\mof=\Ln$ hence $\hx$ gives an minimax estimation of \eqref{eq:kr} state. Last condition in case of stationary matrixes $H(t),C(t)$ means that system~\eqref{eq:kr} is full observable hence this result coincides with well-known theorems of linear systems observability.
\end{document}